\documentclass[a4paper,12pt]{article}

\usepackage[latin1]{inputenc}
\usepackage[T1]{fontenc}
\usepackage{amsthm,amssymb,amsfonts,latexsym,amsmath,amstext,xspace}
\usepackage[english]{babel}
\usepackage[a4paper,portrait]{geometry}

\usepackage{times}

\catcode`\@=11
\newcommand{\languefrancaise}{
\catcode `\?=\active
\def?{\relax\ifhmode\ifdim\lastskip>\z@\unskip\fi
\kern.2em\fi \string?}

\catcode `\;=\active
\def;{\relax\ifhmode\ifdim\lastskip>\z@\unskip\fi
\kern.2em\fi \string;}

\catcode `\:=\active
\def:{\relax\ifhmode\ifdim\lastskip>\z@\unskip\fi
\penalty\@M\ \fi \string:}

\catcode `\!=\active
\def!{\relax\ifhmode\ifdim\lastskip>\z@\unskip\fi
\kern.2em\fi \string!}

\frenchspacing } \catcode`\@=12

\long\def\xcom#1{}

\newcommand{\ladate}{\space\the\day\ \ifcase\month\or janvier\or f\'evrier
\or mars\or avril\or mai\or juin\or juillet\or aout\or septembre
\or octobre\or novembre\or d\'ecembre\fi \ {\oldstyle\the\year}}

\def\alphabet#1{\ifcase#1\or a\or b\or c\or d\or e\or f\or g\or
h\or i\or j\or k\or l\or m\or n\or o\or p\or q\or r\or s\or t\or
u\or v\or w\or x\or y\or z\fi}

\newcommand{\R}{\mathbb{R}}
\newcommand{\PP}{\mathbb{P}}
\newcommand{\E}{\mathbb{E}}

\newcommand{\unsur}[1]{{\frac{1}{#1}}}

\def\valabs#1{{\left\vert {#1} \right \vert}}
\def\esp#1{{{\E}\left [ {#1} \right ]}}
\def\etp#1{{\left ( {#1} \right )}}
\def\etc#1{{\left [ {#1} \right ]}}
\def\crochet#1{{\left < {#1} \right >}}

\def\prob#1{{{\PP}\left ( {#1} \right )}}
\def\norme#1{{\left \Vert #1 \right \Vert}}
\def\ens#1{{\left\{#1\right\}}}

\def\dessus#1#2{\mathord{\mathop{\kern 0pt #2}\limits^#1}}
\newcommand{\egaldef}{{\;\mathrel{\stackrel{\hbox{def}}{=}}\,}}

\newcommand{\bit}{\begin{itemize}}
\newcommand{\eit}{\end{itemize}}
\newcommand{\ben}{\begin{enumerate}}
\newcommand{\een}{\end{enumerate}}

\newcounter{moncompteur}

\newenvironment{myenumerate}%
{\begin{list}{\arabic{moncompteur}. }{\usecounter{moncompteur}%
\setlength{\leftmargin}{0pt}%
\setlength{\labelwidth}{0pt}%
\setlength{\listparindent}{0pt}%
\setlength{\labelsep}{0pt}}}%
{\end{list}}

\def\bmen{\begin{myenumerate}}
\def\emen{\end{myenumerate}}

\newcommand{\undemi}{\frac{1}{2}}

\newcommand{\Frond}{{\mathcal F}}

\newcommand{\Nrond}{{\mathcal N}}

\newcommand{\SK}{Sherrington-Kirkpatrick\xspace}
\newcommand{\gesp}[1]{\mathbf{E}\etc{#1}}
\newcommand{\gprob}[1]{\mathbf{P}\etc{#1}}
\newcommand{\iid}{i.i.d\xspace}

\newtheorem{ethm}{Theorem}

\newtheorem{eprop}[ethm]{Proposition}

\newtheorem{elem}[ethm]{Lemma}

\newtheorem{erem}[ethm]{Remark}

\newcommand{\proofend}{\hfill $\Box{~}$}

\def\as{\rm a.s.\xspace}

\title{Universality in Sherrington-Kirkpatrick's Spin Glass Model}

\author{Philippe {\sc Carmona}\footnote{ P.Carmona: Laboratoire Jean Leray,
    UMR 6629, Université de Nantes, 92208, F-44322, Nantes cedex 03, e-mail:
philippe.carmona@math.univ-nantes.fr},
 Yueyun {\sc Hu}\footnote{Y.Hu: Laboratoire de Probabilités et Modèles
   Aléatoires (CNRS UMR-7599), Université Paris VI, 4 Place Jussieu, F-75252 Paris cedex 05,
e-mail: hu@ccr.jussieu.fr}}

\begin{document}

\maketitle

\hrule
\begin{abstract}
  We show that the limiting free energy in Sherrington-Kirkpatrick's Spin Glass Model does not depend on the environment.
\end{abstract}
\hrule
\bigskip

\section{Introduction}
The physical system is an $N$-spin configuration
$\sigma=(\sigma_1,\ldots,\sigma_N) \in \ens{-1,1}^N$. Each
configuration $\sigma$ is given a Boltzmann weight
$e^{\frac{\beta}{\sqrt{N}} H_N(\sigma) + h \sum_i \sigma_i}$ where
$\beta=\unsur{T}>0$ is the inverse of the temperature, $h$ is the
intensity of the magnetic interaction, $H_N(\sigma)$ is the random
Hamiltonian
$$ H_N(\sigma)=H_N(\sigma,\xi)=  \sum_{1\le i,j\le N} \xi_{ij}\sigma_i \sigma_j\,,$$
and $(\xi_{ij})_{1\le i,j\le N}$ is an \iid family of random
variables, admitting order three moments, which we normalize:
\begin{equation}
\label{ass:xi} \gesp{\xi}=0\,,\quad \gesp{\xi^2}= 1\,,\quad
\gesp{\valabs{\xi}^3}< +\infty\,.
\end{equation}

\medskip
The object of interest is the random Gibbs measure
$$\crochet{f(\sigma)}= \unsur{Z_N} 2^{-N} \sum_\sigma f(\sigma) e^{\frac{\beta}{\sqrt{N}}  H_N(\sigma,\xi) + h \sum_i \sigma_i}\,,$$
and in particular the partition function

$$ Z_N=Z_N(\beta,\xi)= 2^{-N} \sum_\sigma  e^{ \frac{\beta}{\sqrt{N}} H_N(\sigma,\xi) + h \sum_i \sigma_i}\,.$$

We shall denote by $g=(g_{ij})_{1\le i,j\le N}$ an environment of
\iid Gaussian standard random variables ($\Nrond(0,1)$).

\medskip

Recently, F. Guerra and F.L. Toninelli
\cite{1013.82023,1004.82004} gave a rigorous proof, at the
mathematical level, of the convergence of free energy to a
deterministic limit, in a Gaussian environment,
$$  \unsur{N} \log Z_N(\beta,g) \to \alpha_\infty(\beta) \quad \text{\as and in average}.$$

Talagrand~\cite{1014.60050} then proved that one can replace the
Gaussian environment by a Bernoulli environment $\eta_{ij}$,
$\prob{\eta_{ij}=\pm 1}=\undemi$, and obtain the {\em same limit}:
$\alpha_\infty(\beta)$. We shall generalize this result.

\begin{ethm}\label{thm:un}
  Assume the environment $\xi$ satisfies~\eqref{ass:xi}. Then,
$$  \unsur{N} \log Z_N(\beta,\xi) \to \alpha_\infty(\beta) \,\quad \text{\as and in average}.$$
Furthermore, the averages $\alpha_N(\beta,\xi) \egaldef \unsur{N}
\gesp{\log Z_N(\beta,\xi)}$ satisfy
$$ \valabs{\alpha_N(\beta,\xi)-\alpha_N(\beta,g)}\le 9  \gesp{\valabs{\xi}^3}\frac{\beta^3}{\sqrt{N}}\,.$$
\end{ethm}

Therefore the limiting free energy $\alpha_\infty(\beta)$ does not
depend on the environment, hence the {\em Universality} in the
title of this paper : this independence from the particular disorder  was already clear to Sherrington and Kirkpatrick~\cite{sk78} although they had no mathematical proof of this fact (Guerra and Toninelli~\cite{1004.82004} provided a physical proof in the case the environment is symmetric with a finite fourth moment).
 
 Notice eventually that $\alpha_\infty(\beta)$ can be determined in a Gaussian framework where Talagrand~\cite{01981387} recently proved that it is the solution of G. Parisi's variational formula.

\bigskip
The universality property can be mechanically extended to the
ground states, that is the supremum of the families of random
variables:
$$ S_N(\xi)= \sup_{\sigma} \sum_{1\le i,j\le N} \sigma_i \sigma_j \xi_{ij} = \sqrt{N}\lim_{\beta\to +\infty} \unsur{\beta} \log Z_N(\beta,\xi)\,.$$

F. Guerra and F.L. Toninelli \cite{1013.82023,1004.82004} proved
that $N^{-3/2}S_N(g)$ converges as and in average to a
deterministic limit $e_\infty$. Here is the generalization :

\begin{ethm}\label{thm:deux}
  Assume the environment $\xi$ satisfies~\eqref{ass:xi}. Then,
$$ N^{-3/2}S_N(\xi)  \to e_\infty \,\quad \text{\as and in average}.$$
Furthermore, the averages  satisfy, for a universal constant
$C>0$,
$$N^{-3/2} \valabs{\gesp{S_N(\xi)}-\gesp{S_N(g)}}\le C \etp{1+\gesp{\valabs{\xi}^3}} \, N^{-1/6}\,.$$
\end{ethm}

We end this introduction by observing that we do not need the
random variables $\xi_{ij}$ to share the same distribution. They
only need to be independent, to satisfy~\eqref{ass:xi} and such
that $\sup_{ij} \gesp{\valabs{\xi_{ij}}^3} < +\infty$.

\section{Comparison of free energies}
Let us begin with an Integration by parts Lemma.
\newcommand{\nfse}{\norme{F''}_\infty }
\newcommand{\vlxi}{\valabs{\xi}}

\begin{elem}
  \label{lem:ipp}
Let $\xi$ be a real random variable such that $\gesp{\valabs{
\xi}^3}<+\infty$ and $\gesp{\xi}=0$. Let $F:\R\to\R$ be twice
continuously differentiable with $\norme{F''}_\infty=\sup_{x\in\R}
\valabs{F''(x)}<+\infty$. Then
$$\valabs{\gesp{\xi F(\xi)} - \gesp{\xi^2} \gesp{F'(\xi)}} \le \frac{3}{2} \norme{F''}_\infty \gesp{\valabs{\xi}^3}\,.$$
\end{elem}
\begin{proof}
  Observe first, that by Taylor's formula,
  \begin{gather*}
    \valabs{F(\xi)-F(0)-\xi F'(0)} \le \frac{\xi^2}{2} \nfse\,,\\
\valabs{F'(\xi) - F'(0)} \le \vlxi \nfse\,.
  \end{gather*}
Therefore,
\begin{align*}
  \valabs{\gesp{\xi F(\xi)} - \gesp{\xi^2} \gesp{F'(\xi)}} &=
\valabs{\gesp{\xi F(\xi)} - \gesp{\xi^2} \gesp{F'(\xi)} - F(0)\gesp{\xi}}\\
&= \valabs{\gesp{\xi (F(\xi) -F(0) -\xi F'(0))} - \gesp{\xi^2} \gesp{F'(0) -F'(\xi)}} \\
&\le \nfse \etp{\undemi \gesp{\vlxi^3} + \gesp{\vlxi}\gesp{\xi^2}} \\
&\le \nfse \etp{\undemi \gesp{\vlxi^3} +\gesp{\vlxi^3}^{\unsur{3}} \gesp{\vlxi^3}^{\frac{2}{3}}}\\
&\le \frac{3}{2} \norme{F''}_\infty \gesp{\valabs{\xi}^3}\,.
\end{align*}
\end{proof}

\bigskip
In the general framework, $X=(X_1,\ldots,X_d)$ is a random vector
defined on a probability space $(\Omega,\Frond,\PP)$ such that for
any $i$ : $\valabs{X_i}\le 1$. The environment is an \iid family
of random variables $(\xi_1,\ldots,\xi_d)$ defined on
$(\Omega^{(\xi)},\Frond^{(\xi)},\mathbf{P})$, distributed as a
fixed random variable $\xi$ satisfying~\eqref{ass:xi}. The Gibbs
measure, partition function and averaged free energy are thus

\newcommand{\sid}{ \sum_{i=1}^d}
\newcommand{\sidx}{ \sid X_i \xi_i}
\newcommand{\crochetz}[1]{{\crochet{#1}}^{(z)}}

\begin{gather*}
  \crochet{f(X)}=\unsur{Z(\beta,\xi)} \esp{f(X)e^{\beta \sidx}} \\
Z(\beta,\xi)=  \esp{e^{\beta \sidx}}\,,\qquad \alpha(\beta,\xi)=
\gesp{\log Z(\beta,\xi)}\,.
\end{gather*}
Observe that to define $\alpha(\beta,\xi)$ we do not need to
assume exponential moments for the random variable $\xi$, since
$\valabs{\log Z(\beta,\xi)}\le \valabs{\beta} \sid
\valabs{\xi_i}$. We now approximate the derivative of the averaged
free energy:
\begin{elem}
  \label{lem:deriv}
$$ \frac{\partial \alpha(\beta,\xi)}{\partial\beta} = \beta \,\gesp{ \sid (\crochet{X_i^2} -\crochet{X_i}^2)} + 9d\gesp{\valabs{\xi}^3} O(\beta^2)\,,$$
where  $\valabs{O(\beta^2)}\le  \beta^2$.
\end{elem}
\begin{erem}
  In a Gaussian random environment, the integration by parts formula is an exact formula, therefore  the remainder $9d\times\gesp{\valabs{\xi}^3} O(\beta^2)$ vanishes.
\end{erem}
\begin{proof}
  We have
$$ \frac{\partial \alpha(\beta,\xi)}{\partial\beta} =\gesp{\unsur{Z(\beta,\xi)}\esp{ \sidx e^{\beta \sidx}}} = \gesp{\sid \xi_i F_i(\xi_i)}\,,$$
with $F_i(z)=\frac{\esp{X_i e^{\beta X_i z + \psi_i(X)}}}{\esp{
e^{\beta X_i z + \psi_i(X)}}}$ and $\psi_i(X)=\beta \sum_{j\neq i}
X_j \xi_j$ independent of $\xi_i$.

If we define $\crochetz{H}= \frac{\esp{H  e^{\beta X_i z +
\psi_i(X)}}}{\esp{ e^{\beta X_i z + \psi_i(X)}}}$, then
$$ \frac{\partial}{\partial z} \crochetz{H} = \beta\etp{\crochetz{H X_i} - \crochetz{H}\crochetz{X_i}}\,.$$
Hence,
\begin{gather*}
  F_i(z)=\crochetz{X_i}\,,\quad F'_i(z) = \beta\etp{\crochetz{X_i^2} - \etp{\crochetz{X_i}}^2} \\
F''_i(z)= \beta^2\etc{ \crochetz{X_i^3} - 3
\crochetz{X_i^2}\crochetz{X_i} + 2 \etp{\crochetz{X_i}}^3}\,.
\end{gather*}
Since $\valabs{X_i}\le 1$, we have $\norme{F''_i}_\infty \le 6
\beta^2$,  $0\le F'_i(z)\le \beta$ and
$$ F_i(\xi_i) =\crochet{X_i}\,,\quad F'_i(\xi_i) = \beta(\crochet{X_i^2} - \crochet{X_i}^2)\,.$$
We infer from Lemma~\ref{lem:ipp} that since $\gesp{\xi^2}=1$,
$$ \gesp{\crochet{X_i}\xi_i} = \gesp{\xi_i F_i(\xi_i)} = \beta \gesp{ \crochet{X_i^2} - \crochet{X_i}^2} +9\gesp{\valabs{\xi}^3} O(\beta^2)\,,$$
with $\valabs{O(\beta^2)} \le  \beta^2$. Therefore,
$$ \frac{\partial \alpha(\beta,\xi)}{\partial\beta} = \beta \,\gesp{ \sid (\crochet{X_i^2} -\crochet{X_i}^2)} +9 d\gesp{\valabs{\xi}^3}\times  O(\beta^2)\,.$$
\end{proof}

The next step is the comparison of the averaged free energies for
the environments $\xi$ and $g$ (standard normal).
\begin{eprop}
  \label{pro:compa}
For any $\beta \in\R$,
$$ \valabs{\alpha(\beta,\xi) -\alpha(\beta,g)}\le 9 d\gesp{\valabs{\xi}^3} |\beta|^3\,.$$
\end{eprop}
\begin{proof}
  The interpolation technique of F.~Guerra relies on the introduction of a two parameter Hamiltonian:
$$ Z(t,x) = \esp{e^{\sqrt{t} \sid X_i g_i + \sqrt{x} \sid X_i \xi_i}}$$
and averaged free energy $\alpha(t,x)=\gesp{\log Z(t,x)}$ where
the environments $g$ and $\xi$ are assumed to be independent of
each other, $g$ being standard normal. By Lemma~\ref{lem:deriv},
\begin{align*}
  \frac{\partial }{\partial t} \alpha = \undemi \gesp{\sid \crochet{X_i^2} -\crochet{X_i}^2} \\
\frac{\partial }{\partial x} \alpha = \undemi \gesp{\sid
\crochet{X_i^2} -\crochet{X_i}^2} + 9d\gesp{\valabs{\xi}^3}
O(\sqrt{x})\,,
\end{align*}
with $\valabs{O(\sqrt{x})}\le  \sqrt{x}$. We follow the path
$x(s)=t_0-s, 0\le s \le t_0$. Then,
$$\valabs{\frac{d}{ds} \alpha(s,t_0-s)} \le 9 d\gesp{\valabs{\xi}^3} \sqrt{t_0}\,,$$
and thus, integrating on $\etc{0,t_0}$
$$\valabs{\alpha(0,t_0)-\alpha(t_0,0)} \le 9 d\gesp{\valabs{\xi}^3} t_0^{3/2}\,.$$
This is the desired result for $\beta>0$ (take
$\beta=\sqrt{t_0}$). For negative $\beta$, we consider the
environment $-\xi$ instead.
\end{proof}

We shall now estimate the fluctuations of free energy, the
environment is still constructed with i.i.d random variables
$(\xi_1,\ldots,\xi_d)$ satisfying~\eqref{ass:xi}.

\begin{elem}\label{lem:fluctu}
 There exists some universal constant $c>0$ such
that
$$ \gesp{ |\log Z(\beta, \xi) - \alpha(\beta, \xi) |^3}
 \le c  \, \gesp{  |\xi|^3} \, \valabs{\beta}^3\, d^{3/2}.$$
 Consequently, we have $$ {\bf E}\left[ \left|  \sup_{(X_i)} \sum_{i=1}^d X_i \xi_i - {\bf E} \Big( \sup_{(X_i)} \sum_{i=1}^d X_i \xi_i\Big)
 \right|^3 \right] \, \le \,  c\, \gesp{  |\xi|^3} \,   d^{3/2}. $$
\end{elem}

\begin{proof}

We shall use a  martingale decomposition. Let ${\cal F}_k=
\sigma\{\xi_1, ...\xi_k\}, k\ge1,$ be the natural filtration
generated by $(\xi_k)$. Consider the sequence of  martingale
difference

$$ \Delta_j :=\gesp{
\log Z(\beta, \xi) \, \big|\, {\cal F}_j} - \gesp{\log Z(\beta,
\xi) \, \big|\, {\cal F}_{j-1}} \quad 1\le j\le d,$$

with ${\cal F}_0$ the trivial $\sigma$-field. Then

 $$ \log
Z(\beta, \xi) - \alpha(\beta, \xi)  = \sum_{j=1}^d \Delta_j.$$

\noindent  Burkholder's martingale inequality says that   for some
universal constant $c'>0$, $$ {\bf E} \Big| \sum_{j=1}^d \Delta_j
\Big|^3   \le c' \, {\bf E} \left( \sum_{j=1}^d
\Delta_j^2\right)^{3/2}.$$

\newcommand{\qj}{\mathbb{Q}^{(j)}}

\noindent To estimate $\Delta_j$, we define $ Z^{(j)} := \esp{
e^{\beta \sum_{i=1, i\not= j}^d X_i \xi_i}}$ and an auxiliary
random probability  measure $\qj$ by $$ \qj\left( F( X_1,..., X_d)
\right) := {1\over Z^{(j)}}\, \esp{
 F( X_1,..., X_d)  \, e^{\beta \sum_{i=1, i\not= j}^d X_i
\xi_i}}, \qquad \forall \, F(\cdot)\ge0.$$

Then $$ Z(\beta, \xi) = Z^{(j)}\, \qj \left( e^{ \beta X_j\,
\xi_j}\right).$$

Since $Z^{(j)}$ is independent of $\xi_j$, $\log Z^{(j)}$ has the
same conditional expectation with respect to ${\cal F}_j$ as to
${\cal F}_{j-1}$. It follows that $$ \Delta_j = {\bf E} \left(
\log \qj \left( e^{ \beta X_j\, \xi_j}\right) \big|\, {\cal
F}_j\right) - {\bf E} \left( \log \qj \left( e^{ \beta X_j\,
\xi_j}\right) \big|\, {\cal F}_{j-1}\right).$$

\noindent Using the fact that $|X_j|\le 1$, we get  $ \big|  \log
\qj \left( e^{ \beta X_j\, \xi_j}\right) \big|  \le  \beta |\xi_j
|.$  This implies that $$ |\Delta_j |\le \beta \big( |\xi_j| +
{\bf E} |\xi_j|\big).$$

\noindent It follows that \begin{eqnarray*} {\bf E} \Big| \log
Z(\beta, \xi) - \alpha(\beta, \xi) \Big|^3 &\le & c' \, {\bf E}
\left( \sum_{j=1}^d
\Delta_j^2\right)^{3/2} \\
    &\le & c' \, \beta^3\, {\bf E} \left( \sum_{j=1}^d
\big(|\xi_j| + {\bf E} |\xi_j|\big)^2\right)^{3/2} \\
    &\le & c' \, \beta^3\, \sqrt d\,  \sum_{j=1}^d {\bf E} \big(|\xi_j| + {\bf E}
    |\xi_j|\big)^3 \\
    &\le & c  \, {\bf E}
    |\xi|^3 \, \beta^3\, d^{3/2},
\end{eqnarray*}
where we used the convexity of the function $x\to x^{3/2}$ in the
third inequality. Finally, considering  ${1\over\beta}  \log
Z(\beta, \xi)  $ and letting $\beta \to \infty$, we obtain the
second estimate and end the proof.
\end{proof}

\section{Application to \SK's model of spin glass}
Observe that

$$ Z_N(\beta,\xi) = 2^{-N} \sum_\sigma e^{\frac{\beta}{\sqrt{N}} H_N(\sigma,\xi) + h \sum_i \sigma_i} = \esp{e^{\frac{\beta}{\sqrt{N}} H_N(\tau,\xi) + h \sum_i \tau_i}}\,,$$
where $(\tau_i)_{1\le i\le N}$ are \iid with distribution
$\prob{\tau_i =\mp 1}=\undemi$. We get rid of the magnetic field
by introducing tilted laws:
$$ \prob{\tilde{\tau}_i = \pm 1} = \frac{\undemi e^{\pm h}}{\cosh(h)}\,,\text{so that}\quad \esp{f(\tilde{\tau_i})} = \frac{\esp{f(\tau_i)e^{h \tau_i}}}{\esp{e^{h \tau_i}}}\,.$$
With these notations we have
$$ Z_N(\beta,\xi) = \cosh(h)^N \esp{e^{\frac{\beta}{\sqrt{N}} H_N(\tilde{\tau},\xi)}}\,.$$

\subsection*{Convergence of free energy : Theorem 1}

Applying Proposition~\ref{pro:compa} to
$X_{ij}=\tilde{\tau}_i\tilde{\tau}_j$, $\beta\to
\frac{\beta}{\sqrt{N}}$ and $d=N^2$ yields
\begin{align}\label{eq:trois}
 \valabs{\alpha_N(\beta,\xi)-\alpha_N(\beta,g)} &=
\unsur{N} \valabs{\alpha(\frac{\beta}{\sqrt{N}},\xi)-\alpha(\frac{\beta}{\sqrt{N}},g)}\notag \\
&\le  \unsur{N} 9 N^2
\gesp{\valabs{\xi}^3}\etp{\frac{\valabs{\beta}}{\sqrt{N}}}^3 = 9
 \gesp{\valabs{\xi}^3}\frac{\valabs{\beta}^3}{\sqrt{N}}\,.
\end{align}

Furthermore, the fluctuations can be controlled by
Lemma~\ref{lem:fluctu}:

 $$ \gesp{\valabs{ {1\over N} \log Z_N(\beta, \xi) - \alpha_N(\beta, \xi)}^3} \le c\, {\bf E}
    |\xi|^3 \, \valabs{\beta}^3\, N^{-3/2}, $$

\noindent this gives the a.s. convergence by Borel-Cantelli's
Lemma.

\subsection*{Convergence of ground state : Theorem 2}
We have,  restricting the sum to a configuration yielding a
maximum Hamiltonian to get the lower bound,
$$ e^{\frac{\beta}{\sqrt{N}}  S_N(\xi)}\ge Z_N(\beta,\xi)   = 2^{-N} \sum_\sigma e^{\frac{\beta}{\sqrt{N}} H_N(\sigma,\xi) }\ge 2^{-N} e^{\frac{\beta}{\sqrt{N}}  S_N(\xi)}\,.$$

Therefore,

$$\unsur{\sqrt{N}}\gesp{S_N(\xi)} \ge \unsur{\beta}N\alpha_N(\beta,\xi) \ge\unsur{\sqrt{N}} \gesp{S_N(\xi)} -\frac{N \log 2}{\beta}\,.$$

Combining with inequality~\eqref{eq:trois} yields, by taking
$\beta=N^{1/6}$

\begin{align*}
\unsur{N^{3/2}} \valabs{\gesp{S_N(g)}-\gesp{S_N(\xi)}} &\le
\frac{2  \log 2}{\beta } + \unsur{\beta }
\valabs{\alpha_N(\beta,\xi) -\alpha_N(\beta,g)}
\\& \le \frac{2\log 2}{\beta } + C \gesp{\valabs{\xi}^3}\frac{\beta^2}{\sqrt{N}}
\\&\le C'\etp{1+\gesp{\valabs{\xi}^3}} N^{-1/6}\,.
\end{align*}

The almost sure convergence follows in the same way from the
control of fluctuations and Borel-Cantelli's Lemma.

\section{Some Extensions and Generalizations}
\subsection{The $p$-spin model of spin glasses}
The partition function is
$$ Z_N(\beta,\xi) = 2^{-N} \sum_\sigma e^{\frac{\beta}{\sqrt{N^{p-1}}} H_N(\sigma,\xi) + h \sum_i \sigma_i} = \esp{e^{\frac{\beta}{\sqrt{N^{p-1}}} H_N(\tau,\xi) + h \sum_i \tau_i}}\,,$$
where $(\tau_i)_{1\le i\le N}$ are \iid with distribution
$\prob{\tau_i =\mp 1}=\undemi$ (we get rid of the magnetic field
by introducing tilted laws so we assume, without loss in
generality, that $h=0$).

The Hamiltonian is
$$ H_N(\sigma,\xi)= \sum_{1\le i_1,\ldots,i_p\le N} \sigma_{i_1} \ldots \sigma_{i_p} \xi_{i_1 \ldots i_p}$$
wher $\xi_{i_1 \ldots i_p}$ is an iid family of random variables
with common distribution satisfying~\eqref{ass:xi}.

Applying Proposition~\ref{pro:compa} to $X_{i_1\ldots
i_p}=\tilde{\tau}_{i_1}\ldots\tilde{\tau}_{i_p}$, $\beta\to
\frac{\beta}{\sqrt{N^{p-1}}}$ and $d=N^2$ yields
$$ \valabs{\alpha_N(\beta,\xi)-\alpha_N(\beta,g)} \le
9
 \gesp{\valabs{\xi}^3}\frac{\valabs{\beta}^3}{N^{\frac{p-1}{2}}}\,.$$

\subsection{Integration by parts and comparison of free energies}
The more information we get on the random media, the more precise
our comparison results can be. In particular, the more gaussian
the environment looks like, the closer the free energy is to the
gaussian free energy. For example, we shall assume here that  the
random variable $\xi$ satisfies
\begin{equation}
  \label{assbis:xi}
\gesp{\valabs{\xi}^4} < +\infty\,,\quad
\gesp{\xi}=\gesp{\xi^3}=0\,,\quad\gesp{\xi^2}=1\,.
\end{equation}
 A typical variable in this class is the Bernoulli $\gprob{\eta=\pm 1}=\undemi$.

We get the approximate integration by parts formula
\begin{elem}
  Assume that the real random variable $\xi$ satisfies~\eqref{assbis:xi} and that the function $F:\R\to\R$ is of class $C^3$ with bounded third derivative
$\norme{F^{(3)}}_\infty < +\infty$. Then,
$$\valabs{\gesp{\xi F(\xi)} - \gesp{\xi^2}\gesp{F'(\xi)}} \le \norme{F^{(3)}}_\infty \gesp{\xi^4}\,.$$
\end{elem}
\begin{proof}
  This is again Taylor's formula:
  \begin{gather*}
    F(\xi)=F(0) + \xi F'(0) + \undemi \xi^2 F''(0) + O(\valabs{\xi^3} \norme{F^{(3)}}_\infty) \\
F'(\xi)= F'(0) + \xi F''(0) + O(\xi^2 \norme{F^{(3)}}_\infty)\,.
  \end{gather*}
\end{proof}
Repeating, mutatis mutandis, the proof of
Proposition~\ref{pro:compa} we obtain
\begin{eprop}
  \label{pro:compabis}
There exists a constant $C>0$ such that for any environment $\xi$
satisfying~\eqref{assbis:xi}, and for a Gaussian environment $g$,

\begin{equation}
  \label{eq:2}
  \valabs{\alpha(\beta,\xi) -\alpha(\beta,g)}\le C \gesp{\xi^4} d \beta^4\,.
\end{equation}

\end{eprop}

In the framework of Sherrington-Kirkpatrick model of spin glass,
this yields
$$ \valabs{\alpha_N(\beta,\xi) -\alpha_N(\beta,g)}\le C \gesp{\xi^4}\frac{\beta^4}{N}\,.$$

The ground state comparison is now
$$N^{-3/2} \valabs{\gesp{S_N(\xi)}-\gesp{S_N(g)}}\le C \etp{1+\gesp{\valabs{\xi}^4}} \, N^{-1/4}\,.$$

This is  of the same order than  Talagrand's result (Corollary 1.2
of~\cite{1014.60050})  established for  Bernoulli random
variables.

\bigskip


\medskip \par \vskip 1cm

{\em Acknowledgements} We gratefully acknowledge fruitful conversations with Francesco Guerra, held during his stay at Université de Nantes as ``Professeur Invité'' in 2003 and 2004.

\providecommand{\bysame}{\leavevmode\hbox
to3em{\hrulefill}\thinspace}
\providecommand{\MR}{\relax\ifhmode\unskip\space\fi MR }
\providecommand{\MRhref}[2]{%
  \href{http://www.ams.org/mathscinet-getitem?mr=#1}{#2}
} \providecommand{\href}[2]{#2}

\end{document}